\documentclass[11pt]{article}
\usepackage{epsf, amssymb,latexsym,amsmath}

\newcommand{\proof}{\noindent{\bf Proof.\ }}
\newcommand{\qed}{\hfill $\square$ \bigskip}
\newcommand\ch{\hbox{\rm ch}}

\newtheorem{theorem}{\bf Theorem}[section]

\newtheorem{lemma}[theorem]{\bf Lemma}

\newtheorem{observation}[theorem]{\bf Observation}

\def\slika #1{\begin{center}\hskip 0.2mm\epsffile{#1}\end{center}}

\begin{document}

\title{Degenerate and star colorings of graphs on surfaces}

\author{Bojan Mohar\thanks{Supported in part by the ARRS (Slovenia),
   Research Program P1--0297, by an NSERC Discovery Grant and by
   the Canada Research Chair program. On leave from
   Department of Mathematics, IMFM \& FMF, University of Ljubljana,
   Ljubljana, Slovenia}\\
   Department of Mathematics\\Simon Fraser University\\Burnaby, B.C. V5A~1S6\\
   {\tt mohar@sfu.ca}
\and
   Simon \v Spacapan\thanks{Supported in part by the ARRS (Slovenia),
   Research Program P1--0297. The work was written during a post-doctoral 
   visit of the second author at the Simon Fraser University, Canada}\\
   University of Maribor\\FME, Smetanova 17\\2000 Maribor, Slovenia\\
   {\tt simon.spacapan@uni-mb.si}
}

\date{\today}

\maketitle

\begin{abstract}
We study the degenerate, the star and the degenerate star chromatic numbers and their relation to the genus 
of graphs. As a tool we prove the following strengthening of a result of Fertin et al. \cite{fertin}: 
If $G$ is a graph of maximum degree $\Delta$, then $G$ admits a degenerate star coloring using $O(\Delta^{3/2})$ colors. 
We use this result to prove that every graph of genus $g$ admits a degenerate star coloring with $O(g^{3/5})$ colors. 
It is also shown that these results are sharp up to a logarithmic factor. 
\end{abstract}



% ===================================================================

%\newpage

\section{Concepts}

Let $G=(V,E)$ be a graph. An {\em $n$-coloring} of $G$ is a function  
%$$f:V\rightarrow \mathbb{N}\,,$$
$f:V\rightarrow \mathbb{N}$
such that $|f(V)|\leq n$.
We say that $f$ is a {\em proper} coloring if $f(x)\neq f(y)$ for every edge $xy\in E$. 
A {\em color class} $C_i$ of $f$ is the set $f^{-1}(i)$, where $i\in f(V)$.
Two colorings $f$ and $g$ of $G$ are said to be {\em equivalent} if the partitions
of $V$ into color classes of $f$ and $g$ are equal.
Suppose that for each vertex $v\in V(G)$ we assign a \emph{list} $L(v)\subset \mathbb{N}$
of \emph{admissible colors} which can be used to color the vertex $v$.
A \emph{list coloring} of $G$ is a coloring  such that $f(v)\in L(v)$ for each $v\in V$. If for any 
choice of lists $L(v),v\in  V$, such that $|L(v)|\geq k$, 
there exists a proper list coloring of $G$, then we say that $G$ is {\em $k$-choosable}. 
The {\em list chromatic number} of $G$, denoted as $\ch(G)$, is the least $k$, such that $G$ is $k$-choosable.

%We say that $f:V\rightarrow \{1,2,\ldots,n\}$ and $g:V\rightarrow \{1,2,\ldots,n\}$ are {\em equivalent} colorings 
%if they have equal color classes. 

A proper coloring of $G$, such that the union of any two color classes induces a forest, is called an 
{\em acyclic coloring}. The {\em acyclic chromatic number} of $G$, denoted as $\chi_a(G)$, is the least $n$ 
such that $G$ admits an acyclic $n$-coloring.

The notion of a degenerate coloring is a strengthening of 
the notion of an acyclic coloring.
A graph $G$ is {\em k-degenerate} if every 
subgraph of $G$ has a vertex of degree less than $k$. A coloring of a graph 
such that for every $k\ge1$, the union of any $k$ color classes induces 
a $k$-degenerate subgraph is a {\em degenerate coloring}. 
The {\em degenerate chromatic number} of $G$, denoted as $\chi_d(G)$, is the least $n$ 
such that $G$ admits a degenerate $n$-coloring. 

\begin{figure}[htb!]
\epsfxsize=4.0truecm
\slika{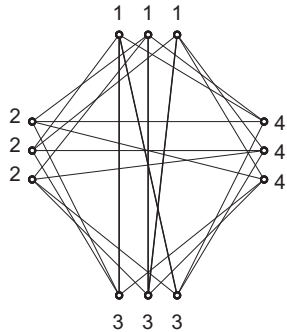}
\caption{An example of a star coloring which is not degenerate}
\label{star-non-degen}
\end{figure}

A proper coloring of $G$, with no two-colored $P_4$ is called a {\em star coloring}. This is 
equivalent to saying that the union of any two color classes induces a star forest, i.e.
a subgraph whose every component is a star $K_{1,t}$ for some $t\ge 0$. The least 
$n$ such that $G$ admits a star coloring with $n$ colors is called the {\em star chromatic number} of $G$, 
denoted as $\chi_s(G)$.

If a coloring is both, degenerate and star, then we speak of 
a {\em degenerate star coloring}. The corresponding 
chromatic number is denoted as $\chi_{sd}$. 

%For all chromatic numbers $\chi_a,\chi_d,\chi_s$ and $\chi_{sd}$ their list versions are defined analogously and 
%we denote them by $\ch_a,\ch_d,\ch_s$ and $\ch_{sd}$. For example, a 
A proper list coloring is an acyclic coloring if the union of any two color classes induces a forest. The {\em 
acyclic list chromatic number} $\ch_a$ is  the least $n$, such that for any assignment of lists of size $n$, there is an  
acyclic list coloring of $G$. The definitions of list versions for all other types of chromatic numbers are 
analogous to their non-list versions and we denote the list versions of chromatic numbers 
by $\ch_a,\ch_d,$ $\ch_s$ and $\ch_{sd}$. 

Clearly, $\chi_a(G)\leq \chi_d(G)\leq \chi_{sd}(G)$ and $\chi_a(G)\leq \chi_s(G)\leq \chi_{sd}(G)$. 
However $\chi_d(G)$ and $\chi_s (G)$ are not comparable. To see this, note that 
the degenerate chromatic number of a tree is two. However, for any tree $T$ which is not a star, 
$\chi_s(T)\geq 3$. In Fig.~\ref{star-non-degen} we give an example of a graph whose star chromatic number is four, 
but has no degenerate four-coloring (since its minimum degree is four). 

It is well known that the list chromatic number of a graph of genus $g$ is $O(g^{1/2})$ (see e.g., \cite{MT}). 
For acyclic colorings, Borodin proved in \cite{borodin} that every planar graph admits
an acyclic 5-coloring and thereby solved a conjecture proposed by Gr\"unbaum \cite{zelenodrevo}. 
Alon et al.~\cite{alon1} determine the (asymptotic) dependence on the acyclic chromatic 
number for graphs of genus $g$, where $g$ is large.
The corresponding bounds for the acyclic list chromatic number have not appeared in
the literature, but the proof in \cite{alon1} can be rather easily adapted to
give the same bounds for the list chromatic version.  

It is also conjectured in \cite{borodin} that every planar graph can be colored with five colors, 
so that the union of any $k$-color classes induces a $k$-degenerate graph for $k=1,\ldots,4$. 
Rautenbach \cite{rauten} proved the existence of degenerate colorings of planar graphs using eighteen colors. 
This result was recently improved to nine colors in \cite{MSZ}. 
%\begin{theorem} \label{9}
%Every planar graph $G$ has ${\rm ch}_d(G)\le 9$.
%\end{theorem}

In  \cite{albert} it was proved that every planar graph admits a star 
coloring with twenty colors and  that the star chromatic number of a graph of genus $g$ is $O(g)$. 

\begin{table}[bht]
    \begin{center}
        \begin{tabular}{|l||c|c|l|}
            \hline
             & Planar & Upper bound & Lower bound \\
             \hline\hline
              $\ch(G)$ & ~~~~5 & $O(g^{1/2})$ & $\Omega(g^{1/2})$  \\
             \hline
             $\ch_a(G)$ & $\le ~7$ & $O(g^{4/7})$ & $\Omega(g^{4/7}/\log(g)^{1/7})$ \\
             \hline
             $\ch_d(G)$ & $\le ~9$ & $O(g^{3/5})$ & $\Omega(g^{4/7}/\log(g)^{1/7})$ \\
             \hline
             $\ch_s(G)$ & $\le 20$ & $O(g^{3/5})$ & $\Omega(g^{3/5}/\log(g)^{1/5})$ \\
             \hline
             $\ch_{sd}(G)$ & -- & $O(g^{3/5})$ & $\Omega(g^{3/5}/\log(g)^{1/5})$ \\
             \hline
        \end{tabular}
    \end{center}

    \caption{Bounds for chromatic numbers in terms of the genus ($g$)}
    \label{table}
\end{table}

The aim of this paper is to establish upper and lower bounds for the degenerate and the star 
list chromatic numbers. We prove that the degenerate star choice number of a graph of 
genus $g$ is $O(g^{3/5})$ thereby improving the bound $O(g)$ given in \cite{albert}. 
We also prove that our bound is sharp up to a logarithmic factor. 
These results in particular solve Problem 3 proposed in \cite[Section~8]{albert}.
The results of this paper and previously known results are collected in Table~\ref{table}.

\section{Probabilistic approach}

In this section we give upper bounds on the degenerate star choice number in terms of 
the maximum degree. The proof uses the probabilistic method in a way similar to
that used in \cite{alon} and \cite{fertin} in the case of acyclic and star colorings. 
It is based on the Lov\'asz Local Lemma but the proof is more complicated than the 
corresponding proofs in \cite{alon} and \cite{fertin}. We refer to \cite{MoRe} for
applications of the probabilistic method to graph colorings.

For $X,Y\subseteq V$ we denote by $E(X,Y)$ the set of edges with one 
endvertex in $X$ and the other in $Y$.

%Let $G=(V_1,E_1)$ and $H=(V_2,E_2)$ be graphs. A function $f:V_1\rightarrow V_2$ is called an {\em isomorphism} 
%of $G$ and $H$ if for every $x,y\in V_1$ we have $f(x)f(y)\in E_2$ if and only if $xy\in E_1$. 
%We say that $G$ and $H$ are {\em isomorphic} if there is an isomorphism from $G$ to $H$, and we write $G\cong H$.

\begin{observation}
\label{obs:1}
Let\/ $G$ be a graph with minimum degree $k\geq 2$ and let $f$ be a proper $k$-coloring of\/ $G$. 
If\/ $S$ is a non-empty subset of a color class $C_i$ of $f$, then there exist a color class 
$C_j$ of $f$, such that $|E(S,C_j)|\geq \frac{k}{k-1}|S|>|S|$.  
\end{observation}

\proof
Each vertex in $S$ has degree at least $k$. Therefore, 
$$\sum_{j\neq i}|E(S,C_j)|\geq k|S|,$$ 
which implies the claimed inequalities. \qed

To prove the main result of this section, we will use the Lov\'asz Local Lemma stated below 
(c.f. \cite{as} or \cite{MoRe}). 

\begin{lemma} \label{local} %(see bla) 
Let $A_1, A_2, \ldots, A_n$ be events in an arbitrary probability space. Let 
$H=(V,E)$ be a graph whose vertices $V=\{1,2,\ldots,n\}$ correspond to the events $A_1, A_2, \ldots, A_n$ 
and whose edge-set satisfies the following: for each $i$, the event $A_i$ is mutually 
independent of the family of events $\{A_j\,|\,ij\notin E\}$. If there exist real 
numbers $0\leq w_i<1$ such that for all $i$
$${\rm Pr}(A_i)\leq w_i\prod_{ij\in E}(1-w_j)$$
then 
$${\rm Pr}\big(\wedge_{i=1}^{n} \bar A_i\big)\geq \prod_{i=1}^n(1-w_i)>0\,,$$
so that with positive probability no event $A_i$ occurs. 
\end{lemma}

Any graph $H$ satisfying the condition stated in the above lemma is called 
a {\em dependency graph} for the events $A_1,\ldots,A_n$. 

The following theorem is the main result of this section. 
It is proved without an attempt to optimize the constant. 

\begin{theorem} \label{main}
For any graph $G$ with maximum degree $\Delta$ there is a degenerate star list coloring 
of\/ $G$ whenever the list of each vertex contains at least 
$\lceil 1000\Delta^{ 3/2}\rceil$ admissible colors. Moreover, a list coloring exists 
such that for every vertex $v$ of degree 
at most $\Delta^{1/2}$, all neighbors of $v$ are colored differently. 
In particular $\chi_{sd}(G)\leq \lceil 1000\Delta^{3/2}\rceil$. 
\end{theorem}

\proof
Let $G$ be a graph with maximum degree $\Delta$ and let $\alpha=\lceil 1000 \Delta^{3/2}\rceil$.   
Suppose that for each vertex $v$ of $G$ a list $L(v)$ of admissible colors is given and that 
$|L(v)|=\alpha$. Consider the uniform probability space  of all list colorings of $G$.  
Then each list coloring of $G$ appears with equal probability.  

We will apply the Lov\'asz Local Lemma to show that in this probability space, a coloring of $G$ 
is a degenerate star list coloring, and has additional properties as stated in the theorem, 
with positive probability. For this purpose, we define events of several types and show that if none of them occurs, 
then the coloring is a degenerate star list coloring with the required properties.

\begin{figure}[htb!]
\epsfxsize=12.0truecm
\slika{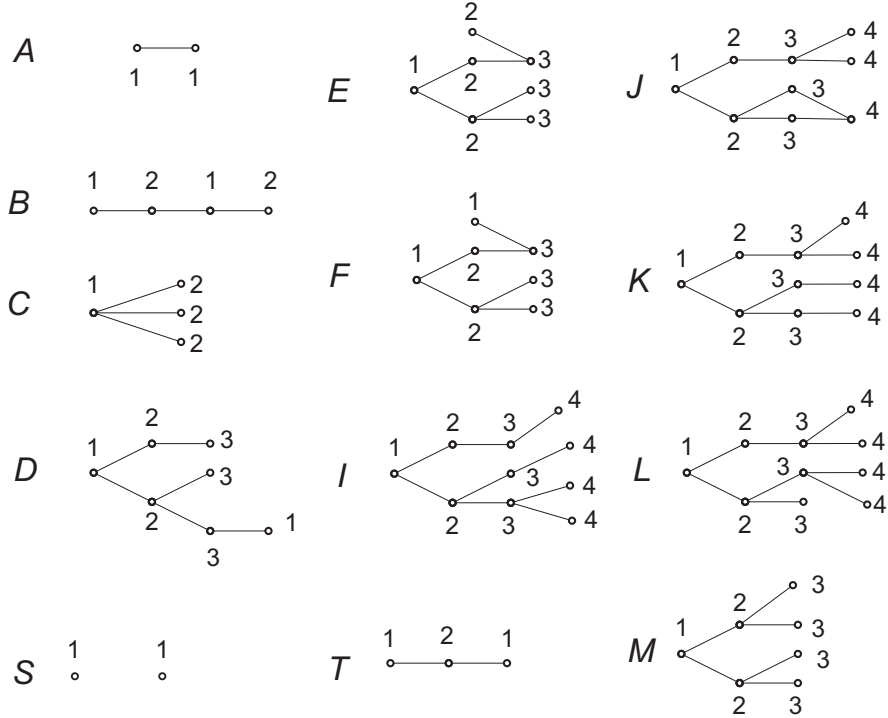}
\caption{The set of graphs $\mathcal{X}=\{A,B,C,D,E,F,I,J,K,L,M,S,T\}$}
\label{forbid}
\end{figure}

A pair of vertices at distance two having at least  
$\Delta^{1/2}$ common neighbors is called a {\em special pair}. 
We next define a family $\mathcal{F}$ of subgraphs of $G$. On Fig.~\ref{forbid} we give a set of graphs $\mathcal{X}$ 
and for each graph $X\in \mathcal{X}$, a coloring $g_X$ (as shown) is given.  
Suppose that $R$ is a subgraph of $G$ isomorphic 
to $X\in \mathcal{X}-\{S,T\}$ 
and let $i: V(X)\rightarrow V(R)$ be an isomorphism of these two graphs. 
If for every pair $a,b\in V(X)$, having the same color under $g_X$, the pair $i(a),i(b)$ is not a special pair 
in $G$, then $R$ is in $\mathcal{F}$. Moreover, $\mathcal{F}$ contains all special pairs
of vertices and all paths $x_1x_2x_3$ such that 
$\deg_G(x_2)\leq \Delta^{1/2}$. 

Let $f$ be a random list coloring of $G$. For $R\in \mathcal{F}$,
denote by $X_R$ the event that the induced coloring $f_{|R}$ on $R$ is equivalent to 
the coloring $g_X$, where $R\in \mathcal{F}, X\in \mathcal{X}$ and $R\cong X$.  
%Additionaly, let $S$ be a disjoint union of two vertices and $T$ a path on three vertices. 
%If $R=\{x,y\}$ is a special pair, let $S_R$ be the event that $f(x)=f(y)$, and if $R=uvw$ is a 
%path on three vertices with $\deg_G(v)\leq 12$, let $T_R$ be the event that $f(u)=f(w)$. 
We refer to an event $X_R$ as an event of type X.  \\

{\sc Claim 0:} If no event $X_R$ occurs for  $X\in \mathcal{X}$, and $R\in \mathcal{F}$,
then the coloring $f$ of $G$ is a degenerate star coloring such that for every vertex 
of degree $\leq\Delta^{1/2}$ all its neighbors are colored by pairwise different colors.\\

{\em Proof:} Let us observe that the definition of the set $\mathcal{F}$ includes
the condition that vertices of the same color do not form a special pair.
However, knowing that events of type S do not occur, we may simply forget about this
condition in the rest of the proof. In particular, if a subgraph $X$ of $G$ is isomorphic
to some $Y\in \mathcal{X}-\{S,T\}$, then its coloring $f_{|X}$ is not equivalent to $g_Y$.

Since no event of type A, B, or T occurs, the coloring $f$ is a star coloring such that any vertex of degree 
at most $\Delta^{1/2}$ has its neighbors colored by pairwise different colors. It remains to prove that the coloring is degenerate. 
Suppose on the contrary, that there is a subgraph $Q$ of $G$ with minimum degree $k$ colored by $k$ colors. 
Since $f$ is a proper coloring, we have $k\geq 2$.
Then there exists a vertex $x\in V(Q)$ of (say) color 1 adjacent to two vertices $y,z$ of color 2 (see Observation \ref{obs:1}).  
Furthermore, there is a color class $P$ of $f_{|Q}$, such that $|E(\{y,z\},P)|\geq 3$. 
Since events of type B do not occur, the color of $P$ is not 1 or 2. Similarly we see that  
$N(y)\cap N(z)\cap P=\emptyset$. As events of type C are excluded, there 
are three vertices $u,v,w\in P$, such that $u,v$ are adjacent to $y$ and $w$ is adjacent to $z$. 
Let $Y=\{u,v,w\}$. Then there is a color class $P'$ such that $|E(Y,P')|\geq 4$. Since  events of type 
B, D, E, and F do not occur, $P'$ is distinct from the color classes 1 and 2. If $|N(Y)\cap P'|\geq 4$ then 
(since B and C do not occur) an event of type I, K, L or M happens, a contradiction. 
If $|N(Y)\cap P'|\leq 3$, then a similar argument shows that either type J or type E event occurs. This contradiction proves the claim. 
\qed
 
Let $H$ be the graph with vertices $X_R$ ($X\in \mathcal{X}, R\in \mathcal{F}, R\cong X)$, 
in which two events $X_{R_1}$ and $Y_{R_2}$ ($X,Y\in \mathcal{X}$) are adjacent if and only 
if $R_1\cap R_2\neq \emptyset$. Since every vertex gets its color independently from others,
$X_R$ is mutually independent of the family of events $Y_{R'}$ 
($Y\in \mathcal{X}, R'\in \mathcal{F}, R'\cong Y, R\cap R'=\emptyset)$.
Therefore, the graph $H$ is a dependency graph for the events $X_R$.

\bigskip
{\sc Claim 1:} Let  $X\in \mathcal{X}, R\in \mathcal{F}$ and  $R\cong X$. 
Then for every $Y\in \mathcal{X}-\{S\}$, the number of events of type Y adjacent to $X_R$ in $H$ 
is at most $100 \Delta^{|Y|-1}$.\\

{\em Proof}:  If  $X_R$ is adjacent to $Y_{R'}$, then there is a vertex $u\in R\cap R'$. 
There are $|X||Y|$ possibilities to choose a one-vertex intersection of $R$ and $R'$. 
Since $G$ has maximum degree $\Delta$ and $R'$ is connected, there are at most $\Delta^{|R'|-1}$ ways to 
choose the other vertices of $R'$. It follows that there are 
at most $|X||Y|\Delta^{|Y|-1}$ graphs isomorphic to $Y$ whose intersection with $R$ is nonempty. 
The result follows from the fact that $|X|\leq 10$ for all $X\in \mathcal{X}$. \qed

{\sc Claim 2:}  Let  $X\in \mathcal{X}, R\in \mathcal{F}$ and $R\cong X$. 
Then the number of events of type J adjacent to $X_R$ in $H$ is at most $90 \Delta^{15/2}$.\\

{\em Proof:} There are $9|R|\leq 90$ possibilities to choose 
a one-vertex intersection of $R$ with a graph $R'$ isomorphic to $J$. Let 
$i:V(R')\rightarrow V(J)$ be an isomorphism and recall that any pair of vertices 
$a,b$ in $R'$ with $g_J(i(a))=g_J(i(b))$ is not a special pair. In particular, a pair of vertices in 
the four-cycle of $R'$ is not a special pair. Therefore there are at most 
$\Delta^{1/2}$ ways to choose a common neighbor of this pair. Since there is a vertex of $R'$ that 
can be chosen in at most $\Delta^{1/2}$ ways, we conclude that there are at most 
$\Delta^{1/2}\Delta^7$ ways to chose the graph $R'$, when a one-vertex intersection of $R$ and $R'$ is fixed. \qed
%Since $|J|=9$ and $|X|\leq 10$, the result follows. \qed 

{\sc Claim 3:} Let  $X\in \mathcal{X}, R\in \mathcal{F}$ and  $R\cong X$. 
Then the number of events of type S adjacent to $X_R$ in $H$ is at most $10\Delta^{3/2}$.\\

{\em Proof:} For a fixed vertex $u\in R$ there are less than $\Delta^2$ induced paths of length 2 
with one endvertex $u$. Each special pair of vertices requires $\Delta^{1/2}$ of these paths. 
Thus there are at most $\Delta^2/\Delta^{1/2}$ vertices $v$, such that $u,v$ is a special pair. 
It follows that the number of special pairs intersecting $R$ is at most $|R|\Delta^{3/2}\leq 10\Delta^{3/2}$. 
\qed 

{\sc Claim 4:} Let  $X\in \mathcal{X}, R\in \mathcal{F}$ and  $R\cong X$. 
Then the number of events of type T adjacent to $X_R$ in $H$ is at most $30\Delta^{3/2}$.\\

{\em Proof:} If $X_R$ and $T_{R'}$ are adjacent in $H$, then $R$ and $R'$ have nonempty intersection. 
There are at most $3|R|$ ways to choose a one-vertex intersection of $R$ and a 
path $x_1x_2x_3$. Since $\deg_G(x_2)\leq \Delta^{1/2}$, the other two vertices of $R'$ 
can be chosen in at most $\Delta^{3/2}$ ways. It follows that there are at most 
$3|R|\Delta^{3/2}\leq 30\Delta^{3/2}$ events of type T adjacent to $X_R$ in $H$. 
\qed

The following table gives upper bounds ${\rm P}(X_R)$ on probabilities of events $X_R$ of different types $X\in \mathcal{X}$:
 \begin{center}
        \begin{tabular}{|c|c|c|c|c|c|}
            \hline
  		Type & A, S, T &  B, C & D, E, F, M & J & I, K, L  \\
            \hline            
            P$(X_R)$ & $\alpha^{-1}$ & $\alpha^{-2}$ & $\alpha^{-4}$ & $\alpha^{-5}$ & $\alpha^{-6}$  \\
            \hline    
    \end{tabular}
    \end{center} 

Let us define the weights for the Local Lemma \ref{local}. 
For each event $X_R$ of type X, let the weight be $w_X=2\, {\rm P}(X_R)$. 
Set $\mathcal{T}=\mathcal{X}-\{J,S,T\}$. To be able to apply Lemma \ref{local} it suffices to show that 
$$
   {\rm P}(X_R) \leq w_X (1-w_J)^{90\Delta^{15/2}}(1-w_S)^{10\Delta^{3/2}}
   (1-w_T)^{30\Delta^{3/2}}\prod_{Y\in \mathcal{T}} \left( 1-  w_Y \right)^{100\Delta^{|Y|-1}}.
$$
To prove this, observe that $1-nx\leq (1-x)^n$, so it suffices to show that 
$$
   \frac 12 \leq (1-90 w_J\Delta^{15/2})(1-10w_S\Delta^{3/2})(1-30 w_T\Delta^{3/2}) 
   \prod_{Y\in \mathcal{T}} ( 1-  100 w_Y \Delta^{|Y|-1}).
$$
But this is easily seen to be true, since the weights for the events are 
$w_J=2\alpha^{-5}\leq 2\cdot 10^{-15}\Delta^{-15/2}$ and 
$w_A=w_S=w_T=2\alpha^{-1}\leq 2\cdot 10^{-3}\Delta^{-3/2}$ and 
$w_Y\leq 2\cdot 10^{-6} \Delta^{1-|Y|}$ for $Y\in \mathcal{X}-\{J,A,S,T\}$. 
%,  we infer that the right side of the above inequality is at least $1/2$, 
%$$ \left(1-90\cdot 2^{-59}\right)\left(1-10\cdot 2^{-11}\right)\left(1-360\cdot 2^{-11}\right)\left(1-100\cdot 2^{-11}\right)(1-100\cdot 2^{-23})^7$$%\geq \frac 12$$
%which completes the proof. 

To conclude, Lemma \ref{local} applies and shows that there exists a coloring $f$ for which
no event $X_R$ ($X\in \mathcal{X}$, $R\in \mathcal{F}$) occurs. 
Finally, Claim 0 shows that $f$ is a coloring whose existence we were to prove.
\qed

\section{Degenerate star colorings for graphs of genus $g$}

An {\em orientation} of a graph $G$ is a function $h:E(G)\rightarrow V(G)$, such that 
for each edge $e=xy\in E(G)$, $h(e)\in \{x,y\}$. We call $h(e)$ the {\em head} of $e$ 
and  the other endvertex of $e$ the {\em tail} of $e$. %, denoted as $t(e)$. 
Let $G$ be a graph and $h$ an orientation of $G$. 
The set of {\em in-neighbors} of $v$ is 
$$N^-(v)=\{u\in V(G)\,|\,uv\in E(G),h(uv)=v\}$$ 
and the set of {\em out-neighbors} of $v$ is
$$N^+(v)=\{u\in V(G)\,|\,uv\in E(G),h(uv)=u\}\,.$$
If $c$ is a coloring of $G$, we define the set of {\em in-colors} $C^-(v)$, and 
{\em out-colors} $C^+(v)$ of $v$ as $C^\pm(v) = \{ c(u) \,|\, u\in N^\pm(v)\}.$
%The set of {\em in-neighbors} of $v$ is the set of all $u\in N(v)$ with $h(uv)=v$ and 
%the set of {\em out-neighbors} of $v$ is the set of all $u\in N(v)$ with $h(uv)=u$.

The following lemma is a simple, yet effective, tool for recognition of star colorings. 
We prove it for the sake of completeness although an equivalent result appears 
in \cite{albert} and some other papers. 

\begin{lemma} \label{star} A proper coloring of $G$ is a star coloring if and only if the edges of G
can be oriented so that for every vertex $v$ %all colors used on the
%in-neighbors of $v$ are pairwise different and 
$$|N^-(v)|=|C^-(v)| { ~~and~~} C^-(v) \cap C^+(v) = \emptyset\,.$$
\end{lemma}

\proof
Suppose that $c$ is a star coloring of $G$. Then the union of any two color classes induces  
a star forest. Orient the edges with one endvertex in color class $A$ and the other 
in color class $B$ so that the tail of each edge is the root of a star induced by $A$ and $B$ 
(the root of a star, which is not a $K_2$, is the vertex of degree at least two, and 
the root of $K_2$ is any vertex of $K_2$). 
It is straightforward that so defined orientation satisfies both conditions from the lemma. 

Conversely, let $c$ be a proper coloring of $G$ and $h$ an orientation of $G$ with properties as 
stated in the lemma. Suppose that $x_1x_2x_3x_4$ is a two-colored path on four vertices and 
(without loss of generality) assume that $h(x_2x_3)=x_3$. Since the in-neighbors of 
$x_3$ are colored by pairwise different colors, we have $h(x_3x_4)\neq x_3$. 
Since $C^-(v) \cap C^+(v) = \emptyset$ we see that $h(x_3x_4)\neq x_4$, a contradiction. 
\qed

The following observation will be used in a recursive construction of degenerate colorings. 

\begin{observation}
\label{miki}
Let $G$ be a graph and let $c$ be a degenerate coloring of a vertex-deleted
subgraph $G-v$. If the neighbors of $v$ are colored by pairwise distinct
colors and we color $v$ by a color which is different from all of those colors,
then the resulting coloring of $G$ is degenerate.
\end{observation}

We are ready for our main result. Its proof is given without intention to optimize
constants. Let us recall that a surface has {\em Euler genus} $g$ if its Euler
characteristic is equal to $2-g$.

\begin{theorem} \label{genus}
Let $G$ be a simple graph embedded on a surface of Euler genus $g$. 
Then $\ch_{sd}(G)\leq \lceil 1000g^{3/5}+100000\rceil.$
\end{theorem}

\proof 
Let $\Sigma$ be a surface of Euler genus $g$, and let $G$ be a graph embedded on $\Sigma$. 
For a vertex $v\in V(G)$, let $L(v)$ be the list of admissible colors. 
We shall prove a stronger statement that $G$ admits a degenerate star list 
coloring from lists of size at least $\alpha:=\lceil 1000g^{3/5}+100000\rceil$, 
such that every vertex of degree $\leq 12$ has its neighbors colored by pairwise distinct colors.  
We will reduce the graph $G$ by using a sequence of edge contractions so that in the resulting 
graph $G_0$ the number of vertices of degree at least 
$\Delta_0:=\lceil \tfrac{1}{4} g^{2/5}+12\rceil $ is at most $\alpha_0:=\lceil 48 g^{3/5}\rceil$. %, where $a$ is a constant to be chosen later. 
We say that a vertex $v$ of $G$ is {\em reducible} if its degree is either at most two,
or its degree is equal to $5-i$ and $v$ is adjacent to a vertex of degree $\leq 9+i$ 
for some $i\in \{0,1,2\}.$

We define a sequence of graphs $G=G_\ell,G_{\ell-1},\ldots, G_0$ as follows. 
We start with $G_\ell = G$, and the precise value of $\ell$ will be determined later. 
Suppose that we have the graph $G_t$. If it has no reducible vertices,
then we stop and adjust the indices so that the current graph is $G_0$ and $\ell$ is 
the number of reductions used to get $G_0$ from $G$. Otherwise, let $v$ be a reducible
vertex of $G_t$. In order to obtain $G_{t-1}$, we perform the following reduction.
If $v$ has degree at most 1, then delete $v$; if $v$ has degree 2, then delete 
$v$ and add an edge between the neighbors of $v$ 
(if not already present). Otherwise, $v$ is of degree $5-i$ where $i\in \{0,1,2\}$. In this case 
contract the edge joining $v$ and a neighbor $u$ of degree at most $9+i$ (and delete possible multiple edges 
that appear after the contraction). The new vertex inherits the list $L(u)$ of admissible colors from $u$.

Let the vertex set of $G_0$ be $\{v_1,v_2,\ldots,v_n\}$, 
where $\deg(v_i)\leq \deg(v_{i+1})$ for $i=1,\ldots,n-1$.\\

\noindent {\sc Claim 0:} The number of vertices in $G_0$ of degree at least $\Delta_0$ is at most $\alpha_0$. \\ %$ag^{3/5}$. 

{\em Proof:} Suppose (reductio ad absurdum) that $\deg(v_k)\geq \Delta_0$, where  $k= n-\alpha_0$. 
Since $G_0$ is a minor of $G$, it has an embedding in the same surface $\Sigma$ of Euler genus $g$ as $G$ has. 
Let $\mathcal{F}$ be the set of faces of this embedding. 
Let $c(v_i)=\deg(v_{i})-6$ for $i=1,\ldots,n$ and for each face $f\in \mathcal F$ let $\bar c(f)=2\deg(f)-6$, where 
$\deg(f)$ is the length of $f$. 
Euler's formula for $\Sigma$ implies (cf., e.g., \cite{MT}) that 
$$
  \sum_{i=1}^n   c(v_i)+\sum_{f\in \mathcal F} \bar c(f)\leq 6g-12\,.
$$ 
Now let $c'(v_i)=\frac 12 \deg(v_{i})$ if $i>k$ and $c'(v_i)=c(v_i)$ otherwise. Then
\begin{eqnarray*}
 \sum_{i=1}^n c'(v_i)+\sum_{f\in \mathcal F} \bar c(f)&=&\sum_{i=1}^n c(v_i)+
 \sum_{i=k+1}^n \left (6-\tfrac 12 \deg(v_i)\right )+\sum_{f\in \mathcal F} \bar c(f) \\
 & \leq& 6g-12+6\cdot 48g^{\frac 35} - 24 g^{\frac 35}\big( \tfrac 14g^\frac25+12\big) < 0\,.
\end{eqnarray*}
Let $c''$ be obtained from $c'$ and $\bar c$ by the following discharging rules,
which preserve the left hand side of the above inequality. 
From each face $f\in \mathcal F$ with $\deg(f)\ge4$ send charge $1$ to each vertex 
of degree 3 lying on the boundary of $f$. If $\deg(f)\ge4$ and three consecutive vertices
on its boundary have degrees $\deg(v_1)\in\{4,5\}$, $\deg(v_2)=11$, and 
$\deg(v_3)\in\{4,5\}$, then $f$ sends charge $1$ to $v_2$ as well. Further, send charge 
$1$ from each vertex $x$ of degree $\geq 11$ to each neighbor $y$ of degree 3 
such that the edge $xy$ is incident with two faces of length $3$,  
and  send $1/2$ to other neighbors of degree 3 and to each neighbor of degree 4 or 5. 
Finally, send $1/5$ from each vertex of degree $10$ to each neighbor of degree 5. 

Our goal is to show that $c''(x)\ge 0$ for every $x\in V(G_0)\cup \mathcal{F}$.
This will imply that
$$
  \sum_{i=1}^n c'(v_i)+\sum_{f\in \mathcal F} \bar c(f) = 
  \sum_{x\in V(G_0)\cup{\mathcal F}} c''(x) \ge 0\, ,
$$
which will in turn contradict the above inequality.

Since $G_0$ has no reducible vertices, a face $f\in \mathcal{F}$ has at most 
$\lfloor \tfrac{1}{2}\deg(f)\rfloor$ incident vertices whose degree is 3, or is 11
(and their neighbors on $f$ have degree 4 or 5).
It sends charge $1$ to them. Therefore,
$$
   c''(f) \ge \bar c(f) - \tfrac{1}{2}\deg(f) \ge 0
$$
whenever $\deg(f)\ge4$. Clearly, $c''(f) = \bar c(f) = 0$ if $\deg(f)=3$.

It is easy to see that each vertex $v$ of degree at least 12 sends charge at most
$\tfrac{1}{2}\deg(v)$, so that $c''(v) \ge c'(v) - \tfrac{1}{2}\deg(v) \ge 0$.
It is also clear that for a vertex of degree ten, $c''(v) \ge c'(v) - \tfrac{1}{5}\cdot 10 > 0$.
If $v$ has degree 11, then $c'(v)=5$. The vertex may send charge $\tfrac{1}{2}$ to
all its neighbors (if they all have degrees 4 or 5). This is the only situation where
its charge $c''(v)$ may become negative. However, in such a case $v$ receives
$1$ from each of its neighboring faces (which are of length at least 4 since there
are no reducible vertices). Therefore, $c''(v)\ge0$ holds in every case.
Vertices of degrees between 6 and 9 have not changed their charge. The reader
will verify that vertices of degrees 3, 4, or 5 have $c''(v)\ge0$ as well.

In conclusion, we have $c''(v)\geq 0$ for every $v\in V(G_0)$ and we have
$c''(f)\geq 0$ for every $f\in \mathcal F$. This yields a contradiction. \qed

Let $k=n-\alpha_0$ be as in the above proof, and let $S=\{v_i\,|\,i>k\}$ be the set of 
{\em special vertices}. 
For each special vertex $v_i$ we choose a color 
$f_{v_i}\in L(v_i)$ so that $f_{v_i}\neq f_{v_j}$ whenever $i\neq j$. 
This is possible since $|L(v_i)|>|S|=\alpha_0$. Define new lists of admissible colors 
for the remaining vertices $v\in V(G)-S$:
$$L'(v)=L(v)-\{f_{v_i}\,|\,i>k\}\,.$$
Note that $|L'(v)|\geq |L(v)|-\alpha_0$ and by Claim 0,  $\deg_{G_0}(v_k)\leq \Delta_0$.
If $\deg_{G_0}(v_k)\leq 144$, then we color the graph induced by vertices $v_1,\ldots,v_k$ 
by a distance-two list coloring with $144^2$ colors, where the distance-two coloring is 
a proper list coloring of $G_0$ such that any two vertices at distance two are colored differently.
If $\deg_{G_0}(v_k) > 144$, then we give the graph induced by vertices 
$v_1,\ldots,v_k$ a degenerate star list coloring with $\lceil 1000 \Delta_0^{3/2}\rceil$ colors, 
where the color for $v_i$ is taken from the list $L'(v_i)$ (see Theorem \ref{main}). 
This is possible since $\max\{144^2, 1000\Delta_0^{3/2}\}\leq |L'(v_i)|$. 
Since $\Delta_0>144$, Theorem \ref{main} assures that each vertex of degree at most 12 has
its neighbors colored by different colors. Then we color 
each special vertex $v_i$ ($i>k$) with the color $f_{v_i}$, which completes the coloring 
of $G_0$. Note that each color used for special vertices has been 
used only once altogether, since it has been deleted from the lists. Therefore 
the coloring of $G_0$ is a degenerate star coloring, 
because the subgraph $G_0-S$ was given a degenerate star coloring. 

Now we extend this coloring to a coloring of $G$. Observe that while contracting an edge $uv$
(where $v$ was a reducible vertex in $G_t$) when going from $G_{t+1}$ to $G_{t}$, 
all vertices preserve their neighbor set except for the common neighbors of $u$ 
and $v$ and for the new vertex (which is again denoted by $u$) 
obtained after the contraction. Observe that $\deg_{G_t}(u)\leq 12$.  

When we go back from $G_t$ to $G_{t+1}$, we add the vertex $v\in V(G_{t+1})-V(G_t)$.
By induction, the coloring of $G_t$ is a degenerate star coloring. Therefore, we can orient 
the edges of $G_t$ so that $|N^-(x)|=|C^-(x)| {\rm ~~and~~} C^-(x) \cap C^+(x) = \emptyset$
for each vertex $x\in V(G_t)$.   

We color the vertex $v$ 
by a color in $L'(v)$ so that the neighbors of each vertex $x\in N_{G_{t+1}}(v)$ with $\deg_{G_{t+1}}(x)\leq 12$ receive 
pairwise different colors and so that the color of $v$ is different from the colors of its neighbors and the colors in   
$C^-(x)$, where $x\in N_{G_{t+1}}(v)$. 
We extend the orientation of $G_t$ to an orientation of $G_{t+1}$ by orienting all edges incident to  $v$ towards $v$. 
Observe that by doing so $|C^-(x)-\{f_{v_i}\,|\, i>k\}|$ remains bounded by $\Delta_0$
for each $x\in V(G_{t+1})-S$.  
Since $|L'(v)| > 4(\Delta_0+1)+12$, there is always an available color. 

We claim that so defined coloring of $G_{t+1}$ is a degenerate star coloring with neighbors of each vertex
of degree at most 12 colored by different colors.
Suppose that $G_t$ was obtained from $G_{t+1}$ by contracting the edge $vu$, where 
$\deg_{G_{t+1}}(v)=5-i$ and $\deg_{G_{t+1}}(u) \le 9+i$ for some $i\in\{0,1,2\}$. 
Then $\deg_{G_{t}}(u)\leq 12$. Therefore $u\in V(G_t)$ has all neighbors colored 
by pairwise different colors and hence also 
$v\in V(G_{t+1})$ has all neighbors colored by pairwise different colors. 
The same is also true if the reduced vertex $v$ is of degree 1 or 2.
It follows from Observation \ref{obs:1} that the coloring of $G_{t+1}$ is degenerate. 
We infer from the definition of the coloring (and orientation) that the set of 
in-colors and the set of out-colors of each vertex are disjoint and, moreover, the in-neighbors of a vertex 
receive pairwise different colors. This shows that the coloring is also a star
coloring and completes the proof. \qed

\section{Lower bounds}

In this section we provide lower bounds which show that the upper bound given 
in Theorem \ref{genus} is asymptotically tight up to a polylogarithmic factor. 
It is interesting to note that the lower bound construction is also based on 
the probabilistic method, although in an essentially different way as in the proof of 
Theorem \ref{main} used to establish the upper bound. 

Let $\mathcal F$ be a family of connected bipartite graphs, each of order at least three.
A proper coloring of $G$ is {\em $\mathcal F$-free} if it contains no two-colored subgraph 
isomorphic to a graph $F\in \mathcal F$. The least $n$ such that $G$ admits an $\mathcal F$-free coloring 
with $n$ colors is denoted by $\chi_{\mathcal F}(G)$. 

Special cases of the following lemma have  appeared in \cite{alon}, \cite{alon1}, and \cite{albert}. 
Here $|G|$ and $||G||$ denote the number of vertices and the number of edges of $G$, respectively.

\begin{lemma} \label{ccc}
Let $\mathcal{F}$ be a family of connected bipartite graphs on at least three vertices and let 
$F\in \mathcal{F}$ be a graph with\/ $\Vert F\Vert$ minimum. Let $G=G_{n,p}$ be the random graph on 
$n$ vertices where each pair of vertices is connected by an edge, randomly and independently, 
with probability $$p=9\Bigl(\frac {\log n}{n}\Bigr)^{1/||F||}\,.$$ 
Then almost surely (that is, with probability tending to one as $n$ tends to infinity), 
$G$ has at most $$9 n^{(2||F||-1)/||F||}(\log n)^{1/||F||} $$%( 5/2 +o(1))$$
edges and $\chi_{\mathcal F}(G)>\tfrac{1}{2|F|}\,n$.
\end{lemma}

\proof It is a standard observation that the number of edges of a random graph with 
edge probability $p$ is almost surely less than $pn^2$ (see, e.g., \cite{as}). 
This proves the claim about $||G||$. 

Let $A\cup B$ be the bipartition of the vertex set of $F$. Set $a=|A|$ and $b=|B|$ and assume that $a\geq b$. 
Suppose that $V_1,\ldots, V_k$ is a partition of the vertex set of $G$ into $k\leq \tfrac{1}{2}n/|F|$ parts. 
Then we delete at most $a-1$ vertices in each $V_i$ so that there is a partition of each $V_i$ into 
sets $U^i_1,\ldots, U^i_{k}$ of size $a$ or $b$ and so that (in all partitions together) 
the number of sets of size $a$ equals (or differs by one) the number of sets of size $b$. 
The number of deleted vertices  will be at most $\tfrac{1}{2}(a-1)n/|F| < n/2$ and therefore there 
are at least $n/2a$ pairwise disjoint sets $U^j_\ell$ 
whose size is $a$ or $b$ and each of them is a subset of a member of the partition $V_1,\ldots, V_k$. 
Since the number of $a$-sets is essentially the same to the number of $b$-sets 
we conclude that there are at least $n/4a$ sets of size $a$ and at least $n/4a$ sets of size $b$.
The probability that the partition $V_1,\ldots, V_k$ is an $\mathcal F$-free coloring 
is at most 
$$
  (1-p^{||F||})^{(n/4a)^2}< {\rm exp}(-(9^{||F||}/16a^2) n\log n)\le n^{-{81 n/64 }}\,.
$$
The last inequality follows from $||F||\geq a\geq 2$. Since there are less than $n^n$ partitions
of the vertex set of an $n$-vertex graph into at most $n/(2|F|)$ classes  
we conclude that the probability that an $\mathcal F$-free coloring with at most $n/(2|F|)$ colors exists 
tends to $0$ as $n$ tends to infinity. 
\qed

\begin{theorem}
For every large enough $g$, there is a graph $G$ embeddable in an orientable (resp. non-orientable) surface of genus $g$, 
such that $\chi_{sd}(G)\geq \chi_s(G)\geq \tfrac{1}{32}g^{3/5}/ (\log g)^{1/5}.$
\end{theorem}

\proof
Let ${\mathcal F}=\{P_4\}$ and let $G$ be a graph of order $n$ with at most 
$9 n^{5/3}(\log n)^{1/3}$ edges and $\chi_s(G)\geq n/8$ (see Lemma \ref{ccc}). 
%Let $S'$ be a surface with minimum genus such that $G$ embeds  
%on $S'$ and denote by $g'$ the genus of $S'$. 
%By Euler's formula we have $n-e+f=2-g'$ and from $n,f\geq 1$  we infer that $g'<9 n^{5/3}(\log n)^{1/3}$. 
It is easy to see that $G$ embeds in an orientable (resp. non-orientable) surface of Euler genus 
$g_0=||G||-1\leq \lceil 9 n^{5/3}(\log n)^{1/3}\rceil=:g$ (in fact, every 2-cell embedding of $G$ satisfies this bound). 
Since $g>n^{5/3}$, it follows that $\log g>\tfrac{5}{3}\log n$. Substituting this to $$g<9 n^{5/3}(\log n)^{1/3}+1$$ 
we conclude that $$n>\tfrac{1}{4}g^{3/5}/(\log g)^{1/5}\,$$ 
(for large enough $g$) and hence the theorem follows from $\chi_s(G)\geq n/8$. 
\qed

%\begin{corollary}
%For every large enough $g$, there is a graph embeddable on a surface of Euler genus $g$, such that 
%$\chi_{sd}(G)\geq (1/24)g^{3/5}/ (\log g)^{1/5}.$
%\end{corollary}

\bigskip
\noindent
{\em Acknowledgement:} We are grateful to Ken-ichi Kawarabayashi for some preliminary
discussions concerning the star chromatic number.

\end{document}